\magnification 1200
\input amssym.def
\input amssym.tex
\parindent = 40 pt
\parskip = 12 pt
\font \heading = cmbx10 at 12 true pt
 at 22 true pt
 at 19 true pt
 at 7 true pt
\def \R{{\bf R}}

\centerline{\heading Resolution of Singularities, Asymptotic Expansions of }
\centerline{\heading Oscillatory Integrals, and Related Phenomena}
\rm
\line{}
\line{}
\centerline{\heading Michael Greenblatt}
\line{}
\centerline{August 6, 2008}
\baselineskip = 12 pt
\font \heading = cmbx10 at 14 true pt
\line{}
\line{}
\noindent{\bf 1. Introduction}

\vfootnote{}{This research was supported in part by NSF grant DMS-0654073}Suppose $f(x)$ is a 
real-analytic function defined in some neighborhood of the origin in $\R^n$. Consider the integral 
$$I_{\lambda} = \int e^{i \lambda f(x)} \phi(x)\,dx \eqno (1.1)$$
Here $\phi(x)$ is a smooth bump function defined on a neighborhood of the origin and $\lambda$ is a real
parameter whose absolute value we assume to be large. If $\nabla f(0) \neq 0$,
then by repeated integrations by parts (see [S] Ch 8 for example), for any $N$ one has an estimate
$|I_{\lambda}| < C_{f,\phi,N}|\lambda|^{-N}$ for appropriate constants $C_{f,\phi,N}$.
In the case where $\nabla f(0) 
= 0$, that is, when $f$ has a critical point at the origin, it can be proven (see [M]) using Hironaka's 
resolution of singularities [H1]-[H2] that if the support of $\phi$ is contained in a sufficiently
small neighborhood of the origin, then as $\lambda \rightarrow \infty$, $I_{\lambda}$ has an 
asymptotic expansion of the form
$$e^{i \lambda f(0)} \sum_{\alpha}\sum_{i = 0}^{n-1} a_{i,\alpha}(\phi) \lambda^{-\alpha}\ln(\lambda)^i
\eqno (1.2)$$
Here the sum in $\alpha$ goes over an increasing arithmetic progression of positive rational numbers,
and the $a_{i, \alpha}$ are distributions with respect to the cutoff $\phi$. We refer to the excellent 
resource [AGV] for further results along these lines. In Theorem 1.2, we will provide 
another proof of the existence of the expansion $(1.2)$ using an elementary resolution of singularities
theorem deriving from [Gr], thereby avoiding the use of Hironaka's theorem or other nonelementary
techniques. Illustrating this elementary method of proving $(1.2)$ and related results while
giving the precise estimates of Theorems 1.1-1.3 can be viewed as the main purpose of this paper. 
In Theorem 1.4, we will also give an elementary proof of
the well-known result of Atiyah [At] and Bernstein-Gelfand [BGe] concerning the meromorphy of 
integrals of $f^z$ for nonnegative real-analytic $f$. It should
be pointed out that there exist other elementary resolution of singularities algorithms which have been 
used for various purposes, notably [BiMi].

The analysis of $(1.1)$ is closely related to the analysis of sublevel set
integrals (see Ch 6-7 of [AGV] for example). Namely, let $f(x)$ be as above and assume $f(0) = 0$. 
We consider the integrals
$$J_t = \int_{\{x: 0 < f(x) < t\}} \phi(x)\,dx \eqno (1.3)$$
Theorem 1.1 below, also proved using an extension of [Gr], will show that for small $t > 0$ we have an asymptotic expansion 
$$J_t \sim \sum_{\alpha}\sum_{i = 0}^{n-1} b_{i,\alpha}(\phi) t^{\alpha}\ln(t)^i
\eqno (1.4)$$ 
Furthermore, each ${\partial^k J_t \over \partial t^k}$ will be seen to have asymptotic
expansion given by termwise differentiation of $(1.4)$. This appears to only have been explicitly
done in the case where $f$ has an isolated zero at the origin [J] [Va], although one should note that
if one is willing to assume Hironaka's results it can be proved without a huge amount of difficulty.
[L] considers some related issues in this subject. 
Once one knows the expansion
$(1.4)$, one can obtain the expansion $(1.2)$ using well-known techniques.
In the integral $(1.1)$, after
factoring out a $e^{i \lambda f(0)}$ if necessary, one may work under the assumption that $f(0) = 0$.
In this situation
one first integrates over a given level set ${f = t}$ and then integrates with respect to $t$. 
Consequently, the integral may be rewritten as 
$$\int_0^{\infty} {\partial J_t \over \partial t} e^{i \lambda t}\gamma(t)\,dt + 
\int_0^{\infty} {\partial \bar{J}_t \over \partial t} e^{-i \lambda t}\gamma(t)\,dt \eqno (1.5)$$
Here $\bar{J}_t$ denotes the analogue of $J_t$ with $f$ replaced by $-f$, and $\gamma(t)$ denotes a
bump function such that $\gamma(f(x)) = 1$ for all $x \in supp(\phi)$.  With a little care, one can
then substitute $(1.4)$ for $f$ and $-f$ respectively into $(1.5)$, integrate term by term, and 
obtain $(1.2)$. We will do this rigorously in the proof of Theorem 1.2 below.

It turns out that there are some natural generalizations of $(1.4)$ which are not any harder to prove
using resolution of singularities than $(1.4)$ itself. For example suppose $A = \{x \in \R^n: g_1(x) 
> 0,...,g_k(x) > 0\}$
where the $g_i$ are real-analytic. Assume further that $0 \in bd(A)$. Then one can obtain an asymptotic
expansion for $\int_{\{x \in A: 0 < f(x) < t\}} \phi(x)\,dx$ in the same fashion that one obtains 
$(1.4)$. Furthermore, if one has several real-analytic functions $f_1(x),...,f_l(x)$ each 
satisfying $f_i(0) = 0$, then one can similarly obtain an asymptotic expansion for 
$\int_{\{x \in A: 0 < f_1(x),\,...\,,0 < f_l(x) < t\}} \phi(x)\,dx$. However, the issue of exhibiting 
an asymptotic expansion for
$\int_{\{x \in A: 0 < f_1(x) < t_1,\,...\,,0 < f_l(x) < t_l\}} \phi(x)\,dx$ in $t_1,...,t_l$ appears to be 
quite a bit harder, and will not be addressed here. Our first
result is the following, to be proven in section 3 using the elementary resolution of singularities
theorem of section 2.

\noindent {\bf Theorem 1.1:} Suppose $f_1(x),...,f_l(x)$ and $g_1(x),...,g_k(x)$ are real-analytic 
functions defined on a 
neighborhood of the origin in $\R^n$, with $f_i(0) = 0$ for all $i$. Let $A = \{x \in 
\R^n: g_1(x) > 0,...,g_k(x) > 0\}$, and assume $0 \in bd(A)$. There is a neighborhood $V$
of the origin such that if $\phi(x)$ is a $C^{\infty}$ function supported in $V$, then 
$J_t = \int_{\{x \in A: 0 < f_1(x) < t,\,...\,,0 < f_l(x) < t\}} \phi(x)\,dx$ has an asymptotic 
expansion given by
$$J_t = \sum_{\alpha \leq a}\sum_{i = 0}^{n-1} b_{i,\alpha}(\phi) t^{\alpha}\ln(t)^i + E_a(t)
\eqno (1.6)$$
Here the $\alpha$ range over an arithmetic progression of positive rational numbers depending on
the $f_i$ and the $g_i$. Let $Z = 
\{x \in cl(A): f_i(x) = 0$ for some $i.\}$ There are $M > 0$ and $A_{\alpha} > 
0$ (depending on the $f_i$ and $g_i$) such that each $b_{i, \alpha}$ is a distribution with respect to $\phi$, supported on $Z$ and 
satisfying
$$|b_{i, \alpha}(\phi)| \leq A_{\alpha} \sup_{|\beta| \leq M\alpha}\, \sup_{x \in Z} |\partial^{\beta}
\phi(x)|\eqno (1.7a)$$
The error term $E_a(t)$ is such that there is $\epsilon > 0$ and $C_a > 0$ 
such that if $0 \leq m \leq a$, then
$$|{d^m  \over dt^m}E_a(t)| < C_a \sup_{|\beta| \leq Ma}\, \sup_{x } |\partial^{\beta}
\phi(x)|\,t^{a + \epsilon - m} \eqno (1.7b)$$

\noindent The theorem of this paper regarding oscillatory integrals is as follows:

\noindent {\bf Theorem 1.2:} Suppose $f(x)$ is a real-analytic function defined on a 
neighborhood of the origin in $\R^n$ with $f(0) = 0$, and let $A$ and $V$ be as in Theorem 1.1. Then as 
$\lambda \rightarrow \infty$, 
$I_{\lambda} = \int_A e^{i \lambda f(x)} \phi(x)\,dx$ has an asymptotic expansion 
$$\sum_{\alpha \leq a}\sum_{i = 0}^{n-1} a_{i,\alpha}(\phi) \lambda^{-\alpha}\ln(\lambda)^i
+ E_a'(\lambda)\eqno (1.8)$$
The $\alpha$ range over an arithmetic progression of positive rational numbers depending on
$f$ and the $g_j$.
Let $Z = \{x \in cl(A): f(x) = 0\}$. There are $M' > 0$ and $A_{\alpha}' > 
0$ such that each $a_{i, \alpha}$ is a distribution with respect to $\phi$, supported on $Z$ and 
satisfying
$$|a_{i, \alpha}(\phi)| \leq A_{\alpha}' \sup_{|\beta| \leq M'\alpha}\, \sup_{x \in Z} |\partial^{\beta}
\phi(x)|\eqno (1.9a)$$ 
The error term $E_a'(t)$ is such that there is $\epsilon' > 0$ and $C_{m,a}' > 0$ 
such that for any $m$ one has 
$$|{d^m  \over d\lambda^m}E_a'(\lambda)| < C_{m,a}' \sup_{|\beta| \leq M'(a + m)}\, \sup_{x}\, |\partial^{\beta}
\phi(x)|{\lambda}^{-a - \epsilon' - m} \eqno (1.9b)$$
\noindent {\bf Proof:} As in $(1.5)$ we write
$$I_{\lambda} = \int_0^{\infty} {\partial J_t \over \partial t} e^{i \lambda t}\gamma(t)\,dt + 
\int_0^{\infty} {\partial \bar{J}_t \over \partial t} e^{-i \lambda t}\gamma(t)\,dt \eqno (1.10)$$
As before $\bar{J}_t$ denotes the analogue of $J_t$ with $f$ replaced by $-f$ and $\gamma(t)$ has compact
support and is equal to 1 on a neighborhood of the origin. The two terms of 
$(1.10)$ are done in the same fashion, so we focus our attention on the first term. 
By Theorem 1.1, for a given $k$ we can write 
$${d J_t \over dt} = \sum_{\alpha \leq k-1}\sum_{i = 0}^{n-1} B_{i,\alpha}
(\phi) t^{\alpha}\ln(t)^i + {d E_k \over dt} (t)$$
Here the $B_{i,\alpha}$ are obtained from the $b_{i,\alpha}$ by performing the appropriate term-by-term
differentiation. In order to be able to differentiate the expansion $(1.8)$ up to 
$\lambda^{-a}$ a total of $m$ times, $m \geq 0$, we need to insert the above expansion for $k = 
a + m + 2$ into the term of $(1.10)$. We get 
$$I_{\lambda} = \sum_{\alpha \leq a + m + 1}\sum_{i = 0}^{n-1} B_{i,\alpha}(\phi) \int_0^{\infty} 
t^{\alpha}\ln(t)^i
e^{i \lambda t}\gamma(t)\,dt + \int_0^{\infty} {d E_{a + m + 2} \over dt}(t)e^{i \lambda t}
\gamma(t)\,dt \eqno (1.11)$$
It is well-known (see [F]) that for any $l > 0$ one has 
$$\int_0^{\infty} e^{i \lambda t}t^{\alpha}\ln(t)^m \gamma(t)\,dt = {\partial^m \over \partial \alpha^m} 
{\Gamma(\alpha+1) \over (-i\lambda)^{\alpha + 1}} + O(\lambda^{-l}) \eqno (1.12)$$
As a result, for any $l$, $(1.11)$ becomes
$$I_{\lambda} = \sum_{\alpha \leq a + m + 1}\sum_{i = 0}^{n-1} B_{i,\alpha}(\phi) {\partial^i \over 
\partial \alpha^i} {\Gamma(\alpha+1) \over (-i\lambda)^{\alpha + 1}}
 + \int_0^{\infty} {d E_{a + m + 2} \over dt}(t)e^{i \lambda t}\gamma(t)\,dt + O(\lambda^{-l}) \eqno (1.13)$$
Equation $(1.13)$ will give the desired expression for the  $m$th derivative of $I_{\lambda}$'s 
expansion up to order $\lambda^{-a}$.
For in $(1.12)$ the $O(\lambda^{-l})$ behaves as needed under differentiation (again see [F]).
As for the ${d E_{a + m + 2} \over dt}$ term in $(1.13)$, one can differentiate $\int_0^{\infty} 
{d E_{a + m + 2} \over dt}(t)e^{i \lambda t}\gamma(t)\,dt$ under the integral sign, obtaining a 
power of $it$
for each of the $m$ $\lambda$-derivatives taken. Then one does $a+ m + 1$ integrations by parts in $t$, 
integrating
the $e^{i \lambda t}$ factor and differentiating the rest. Equation $(1.7b)$ ensures the left
endpoint terms disappear. We obtain
$$|{d^m \over d \lambda^m} \int_0^{\infty} {d E_{a + m + 2}\over d t}(t)
e^{i \lambda t}\gamma(t)\,dt| \leq C_{a,m}\,\lambda^{-a - m - 1}\sup_t |{d^{a + m + 2}E_{a + m + 2} \over dt^{a + m + 2}}(t)|$$
$$\leq C_{a,m}'\sup_{|\beta| \leq M'(a + m)}\, \sup_{x}\, 
|\partial^{\beta} \phi(x)| \lambda^{-a - m - 1} \eqno (1.14)$$
The last inequality follows from the ${d^{a + m + 2}E_{a + m + 2} \over dt^{a + m + 2}}(t)$ case 
of $(1.7b)$.
Inserting this back into $(1.13)$ gives the desired result. This completes the proof of Theorem 1.2.

Next, we focus on the situation where $f \geq 0$ on $A$. Where $\tau$ is now a large positive 
parameter, we 
consider the Laplace Transform-like object defined by
$$L_{\tau} = \int_A e^{- \tau f(x)} \phi(x)\,dx \eqno (1.15)$$
\noindent We have the following theorem regarding $L_{\tau}$.

\noindent {\bf Theorem 1.3:} Suppose $f(x)$ is a real-analytic function defined on a 
neighborhood of the origin in $\R^n$ with $f(0) = 0$, and let $A$, $V$, and $Z$ be as in Theorem 1.1. 
Suppose on
a sufficiently small neighborhood $U$ of the origin we have that $f(x) \geq 0$ on $U \cap A$. Then if
the support of $\phi$ is contained in $V \cap U$, as $\tau 
\rightarrow \infty$, $L_{\tau}$ has an asymptotic expansion 
$$L_{\tau} = \sum_{\alpha \leq a}\sum_{i = 0}^{n-1} d_{i,\alpha}(\phi) \tau^{-\alpha}\ln(\tau)^i
+ E_a''(\tau)\eqno (1.16)$$
The $\alpha$ range over an arithmetic progression of positive rational numbers depending on
$f$ and the $g_j$. There are $M'' > 0$ and $A_{\alpha}'' > 
0$ such that each $d_{i, \alpha}$ is a distribution with respect to $\phi$, supported on $Z$ and 
satisfying
$$|d_{i, \alpha}(\phi)| \leq A_{\alpha}'' \sup_{|\beta| \leq M''\alpha}\, \sup_{x \in Z} |\partial^{\beta}
\phi(x)| \eqno (1.17)$$ 
The error term $E_a''(\tau)$ is such that there is $\epsilon'' > 0$ and $C_{m,a}'' > 0$ 
such that for any $m$ one has 
$$|{d^m  \over d\tau^m}E_a''(\tau)| < C_{m,a}''\sup_{|\beta| \leq M''(a + m)}\, \sup_x |\partial^{\beta}
\phi(x)|\,{\tau}^{-a - \epsilon'' - m} \eqno (1.18)$$

\noindent {\bf Proof:} The proof is essentially a repeat of that of Theorem 1.2. Namely, we have
$$L_{\tau} = \int_0^{\infty} {\partial J_t \over \partial t} e^{-\tau t}\gamma(t)\,dt \eqno (1.19)$$
Because of the $f \geq 0$ condition, there is only one term instead of two this time. 
Instead of using $(1.12)$, here we use
$$\int_0^{\infty} e^{- \tau t}t^{\alpha}\ln(t)^i \gamma(t)\,dt = {\partial^i \over \partial \alpha^i} 
{\Gamma(\alpha+1) \over \tau^{\alpha + 1}} + O(\tau^l) \eqno (1.20)$$
Otherwise the proof is identical to that of Theorem 1.2 so we omit the details. 
Incidentally, equation $(1.20)$ is somewhat easier to prove than $(1.12)$. When $i = 0$ one can write 
$(1.20)$ as as 
$$\int_0^{\infty} e^{- \tau t}t^{\alpha}\,dt - \int_0^{\infty} e^{- \tau t}t^{\alpha}(1 - \gamma(t))
\,dt \eqno (1.21)$$
The first term is exactly ${\Gamma(\alpha + 1) \over \tau^{\alpha + 1}}$. On the other hand, via 
repeated integrations by parts, the 
second term and its $\tau$ derivatives are seen to be $O(\tau^l)$ for any $l$. 
The $i > 0$ case follows  from  differentiating $(1.21)$ under the integral with respect to 
$\alpha$.

Theorem 1.1 also gives as a relatively straightforward consequence the following theorem of Atiyah [At]; 
a similar result is due to Bernstein-Gelfand [BGe]:

\noindent {\bf Theorem 1.4:} Let $f$, $A$, $U$, and $V$ be as in Theorem 1.3. 
Define $F_{\phi}(z) = 
\int_A f(x)^z \phi(x)\,dx$. Then if the support of $\phi$ is contained in $V \cap U$, the function 
$F_{\phi}(z)$, initially defined as a holomorphic function of $z$ on 
$Re(z) > 0$, extends to a meromorphic function on all of ${\bf C}$. The poles of this extension
are located on an arithmetic progression of negative rational numbers depending on
$f$ and the $g_j$, and each pole is of order at most the dimension $n$. 

\noindent {\bf Proof:} Replacing $f$ by $cf$ for an appropriate $c > 0$ if necessary, we may assume
that $|f| \leq 1$ on $supp(\phi)$. Analogous to $(1.19)$, we have
$$F_{\phi}(z) = \int_0^1 t^{z} {d J_t \over dt}\,dt \eqno (1.22)$$
Inserting $(1.11)$ into $(1.22)$, one obtains 
$$F_{\phi}(z) = \sum_{\alpha \leq a}\sum_{i = 0}^{n-1} B_{i,\alpha}(\phi) \int_0^1 t^{\alpha + z}
\ln(t)^i\,dt + \int_0^1 {d E_a \over dt}(t) t^z\,dt \eqno (1.23)$$
When $Re(z) > -1$, each of the integrations of the sum in $(1.23)$ can be performed directly, and we 
obtain
$$F_{\phi}(z) = \sum_{\alpha \leq a}\sum_{i = 1}^n B_{i,\alpha}'(\phi) (z + \alpha + 1)^{-i} + 
\int_0^1 {d E_a \over dt}(t)\,t^z\,dt \eqno (1.24)$$
Note that the sum in $(1.24)$ automatically extends to a meromorphic function on ${\bf C}$ with poles of
order at most $n$. As for the 
error term, one can rewrite it as $\int_0^1 F_a(t)t^{z + a -1 +\epsilon}\,dt$, where in view of $(1.7b)$ 
$F_a(t)$ is bounded. Hence on $Re(z) > -a + 1 - \epsilon$, the error term is an average of locally 
uniformly bounded analytic functions. Hence it is itself an analytic function on $Re(z) > 
 -a + 1 - \epsilon$. Since $a$ can be made arbitrarily large, the theorem follows.
 
The remainder of the paper is organized as follows. In section 2, a version of the elementary resolution
of singularities algorithm of [Gr] will be developed that will be appropriate for proving the type of 
theorems of this paper. It may also be useful for other purposes as well. Section 3 will be devoted to 
proving Theorem 1.1 using the algorithm of section 2.

\noindent {\bf 2. A resolution of singularities theorem}

\noindent We now give some terminology from [Gr]:

\noindent {\bf Definition:} We say that a function $g: A \subset \R^n \rightarrow \R^n$ a 
{\it quasitranslation} if there is a real analytic function $r(x)$ of $n-1$ variables such that 
$g(x) = (g_1(x),...,g_n(x))$, where for some $j$ we have $g_j(x) = x_j - r(x_1,...
x_{j-1},x_{j+1},...,x_n)$ and where $g_i(x) = x_i$ for all $i \neq j$.
 In other words $g(x)$ is a translation in the $x_j$ variable when the others are fixed.

\noindent {\bf Definition:} We call a function $m:A \subset \R^n \rightarrow \R^n$ an 
{\it invertible monomial map} if there are nonnegative
integers $\{\alpha_{ij}\}_{i,j=1}^n$ such that the matrix $(\alpha_{ij})$ is invertible 
and $m(x) = (m_1(x),...,m_n(x))$ where $m_i(x) = x_1^{\alpha_{i1}}....x_n^{\alpha_{in}}$.
The matrix $(\alpha_{ij})$ being invertible ensures that $h$ is a bijection on $\{x : x_l > 0
\hbox { for all } l\}$. 

\noindent {\bf Definition:} Let $E = \{x : x_i > 0$ for all $i\}$. If $h(x)$ is a bounded, nonnegative, 
compactly supported function on $E$, we say $h(x)$ is a {\it quasibump function} if $h(x)$ is of the
following form:
$$h(x) = a(x) \prod_{l=1}^j b_l (c_l(x) {p_l(x) \over q_l(x)}) \eqno (2.1)$$
Here $p_l(x), q_l(x)$ are monomials, $a(x) \in C^{\infty}(cl(E))$, the $c_l(x)$ are nonvanishing 
real-analytic functions defined on a neighborhood of $supp(h)$, and $b_l(x)$ are nonnegative functions 
in $C^{\infty}(\R)$ such that there are constants $c_1 > c_0 > 0$ with each $b_l(x) = 1$ for $x < c_0$ and 
$b_l(x) = 0$ for $x > c_1$. 

\noindent The main theorem from [Gr] is as follows:

\noindent {\bf Main Theorem of [Gr]:} Let $f(x)$ be a real-analytic function defined in a neighborhood
of the origin in $\R^n$. Then there is a neighborhood $U$ of the origin such that if $\phi(x) 
\in C_c^{\infty}(U)$ is nonnegative with $\phi(0) > 0$, then 
$\phi(x)$ can be written (up to a set of measure zero) as a finite sum $\sum_i \phi_i(x)$ of
nonnegative functions such that for all $i$, $0 \in supp(\phi_i)$ and $supp(\phi_i)$ is a subset of 
one of the $2^n$ closed quadrants defined by the hyperplanes $\{x_m = 0\}$. The following 
properties hold.

\noindent (1) For each $i$ there are bounded open sets $D_i^0$,...,$D_i^{k_i}$, and maps 
$g_i^1$,..., $g_i^{k_i}$, each a reflection, translation, invertible monomial map, or quasitranslation, 
such that $D_i^0 = \{x: \phi_i(x) > 0\}$ and such that each $g_i^j$ is a real-analytic 
diffeomorphism from $D_i^j$ to $D_i^{j-1}$. The function $g_i^j$ extends to 
a neighborhood $N_i^j$ of $cl(D_i^j)$ with $g_i^j(N_i^j) \subset N_i^{j-1}$ for $j > 1$ and
$g_i^1(N_i^1) \subset U$. \parskip=3pt

\noindent (2) Let $E = \{x : x_i > 0$ for all $i\}$ and $\Psi_i = g_i^1 \circ .... \circ g_i^{k_i}$.
Then $D_i^{k_i} \subset E $, and there is a quasibump function
$\Phi_i$ such that $ \chi_{D_i^{k_i}}(x)(\phi_i \circ \Psi_i(x)) =  \Phi_i(x)$.

\noindent (3) $0 \in N_i^{k_i}$ with $\Psi_i(0) = 0$.

\noindent (4) On $N_i^{k_i}$, the functions $f \circ \Psi_i$, $det(\Psi_i)$, and each $j$th component 
function $(\Psi_i)_j$ is of the form $c(x)m(x)$, where $m(x)$ is a  monomial and $c(x)$ is
nonvanishing.

\parskip = 12pt
The corollary to the main theorem of [Gr] says that one can resolve several functions 
simultaneously in such a way that the resolution satisfies the conclusions of the main theorem. However, 
these theorems are not precisely what is needed for the arguments of this paper, because here a quasibump
function is not the appropriate form for the function $\Phi_i(x)$ of part (2). Instead we will 
need the following:

\noindent {\bf Theorem 2.1:} If in the main theorem of [Gr] and its corollary, if one replaces 
the condition in part (2) that $\Phi_i(x) = 
\chi_{D_i^{k_i}}(x)(\phi_i \circ \Psi_i(x))$ is a quasibump function with the condition 
that for some rectangle $R_i = (0,a_1^i) \times ... \times (0,a_n^i)$ the function $\Phi_i(x)$ is of 
the form $\chi_{R_i}(x) \gamma(x)$, where
$\gamma(x)$ is a $C^{\infty}$ function on $cl(R_i)$, then the rest of the main theorem and its 
corollary respectively still holds.

\noindent One way of looking at Theorem 2.1 is that in the blown-up coordinates, one can replace the issues coming
from the singularities of the quasibump function by the issues coming from the jumps in the 
characteristic function of $R_i$. These latter issues will turn out to cause no problems in the 
analysis of the integral quantities of this paper. Theorem 2.1 will be a consequence of the following,
which we will prove later in this section.

\noindent {\bf Theorem 2.2:} Let $\Phi(x)$ be a
quasibump function, and let $p_l(x)$ and $q_l(x)$ be as in the definition $(2.1)$ of quasibump 
function applied to $\Phi(x)$. Then there is a $J$ such that for $1 \leq j \leq J$ and 
$1 \leq k \leq n$ there are invertible monomial maps $g_j(x)$, 
quasibump functions $Q_j(x)$ of the form $\prod_{l=1}^m \alpha ({p_{jl}(x) \over q_{jl}(x)})$ and sets 
$F_j = \{x \in E: {r_{jk}(x) \over s_{jk}(x)} < 1$ for 
$1 \leq k \leq n\}$ where $r_{jk}(x)$ and $s_{jk}(x)$ are monomials, such that the following hold.

\noindent (1) Up to a set of measure zero one has a decomposition 
$$\Phi(x) = \sum_{j = 1}^{J} \Phi(x)Q_{j}(x) \chi_{F_{j}}(x) \eqno (2.2)$$

\noindent (2) Each ${p_{jl}(g_j(x)) \over q_{jl}(g_j(x))}$ and 
each ${p_l(g_j(x)) \over q_l(g_j(x))}$ is a monomial.

\noindent (3) For each $j$ there is some rectangle $R_j = (0,a_1^j) \times ... \times (0,a_n^j)$
such that $\chi_{F_j}(g_j(x)) = \chi_{R_j}(x)$.

\noindent {\bf Comment:} Note that $\Phi(g(x)) = a(g_j(x)) \prod_{l=1}^L b_l (c_l(g_j(x)) {p_l(g_j(x)) \over
q_l(g_j(x))})$, $a$ and $c_l$ smooth, and that $Q_j(g_j(x)) = \prod_{l=1}^m \alpha ({p_{jl}(g_j(x)) 
\over q_{jl}(g_j(x))})$. Hence (2) implies that $\Phi(g_j(x))$ and each $Q_j(g_j(x))$ are smooth. Thus
by (3), for each $j$ there is and a smooth $\gamma_j(x)$ on $R_j$ with 
$$\Phi(g_j(x))Q_j(g_j(x)) \chi_{F_j}(g_j(x)) = \chi_{R_j}(x) \gamma_j(x) \eqno(2.3)$$

\noindent {\bf Proof that Theorem 2.2 implies Theorem 2.1:}

\noindent Suppose Theorem 2.2 is known to hold, and let $f(x)$ and $U$ be as in the main theorem (or 
its corollary).
Suppose $\phi(x)$ is a bump function defined in $U$, and let $\phi = \sum_i \phi_i$ be the 
decomposition of the main theorem (or its corollary). Let  $\Phi_i(x)$ be the  quasibump function and $\Psi_i$ the
composition of coordinate changes associated to $\phi_i$.  Let $Q_{ij}(x)$, $F_{ij}$, and $g_{ij}$
be as given by Theorem 2.2 applied 
to $\Phi_i(x)$. I claim that the decomposition $\phi(x) = 
\sum_{ij} \phi_i(x) Q_{ij}(\Psi_i^{-1}(x))\chi_{F_{ij}}(\Psi_i^{-1}(x))$ satisfies the conditions of
the main theorem of [Gr], with coordinate changes $g_i^1,....,g_i^{k_i}, g_{ij}$, if one specifies
the domains $D_i$ and $N_i$ as follows. At the $k_i$th level, the domain $D_{ij}^{k_i}$ is  defined
to be $D_i^{k_i} \cap F_{ij} \cap \{x : Q_{ij}(x) > 0\}$. For $k < k_i$, the domains are sucessively
defined by $D_{ij}^{k-1} = g_i^k(D_{ij}^k)$. At the final $k_i + 1$th level, one puts $D_{ij}^{k_i + 1}
= g_{ij}^{-1}D_{ij}^{k_i}$. For 
$j \leq k_i$ one can define $N_{ij}^{k_i}$ to just be $N_i^j$, and then $N_{ij}^{k_i+1}$ to be 
$g_{ij}^{-1} N_{ij}^{k_i}$. 

With the above definitions, part (1) of the main theorem is readily seen to hold. As for (2), 
since $ \chi_{D_i^{k_i}}(x)(\phi_i \circ \Psi_i(x)) =  \Phi_i(x)$, we have
$$\chi_{D_i^{k_i}}(x)(\phi_i \circ \Psi_i(x))Q_{ij}(x)\chi_{F_{ij}}(x) =  \Phi_i(x)Q_{ij}(x)
\chi_{F_{ij}}(x) \eqno (2.4)$$
Therefore,
$$\chi_{D_i^{k_i} \cap \{x: Q_{ij}(x) > 0\} \cap F_{ij}}(x)(\phi_i \circ \Psi_i(x))Q_{ij}(x)
\chi_{F_{ij}}(x) =  \Phi_i(x)Q_{ij}(x)\chi_{F_{ij}}(x)$$
Equvialently,
$$\chi_{D_{ij}^{k_i}}(x) (\phi_i \circ \Psi_i(x))Q_{ij}(x)\chi_{F_{ij}}(x) =  
\Phi_i(x)Q_{ij}(x)\chi_{F_{ij}}(x)$$
Composing with $g_{ij}$, we have
$$\chi_{D_{ij}^{k_i+1}}(x) (\phi_i \circ \Psi_i \circ g_{ij}(x))Q_{ij}(g_{ij}(x))\chi_{F_{ij}}(g_{ij}(x)) =  
\Phi_i(g_{ij}(x))Q_{ij}(g_{ij}(x))\chi_{F_{ij}}(g_{ij}(x)) \eqno (2.5)$$
If one lets $\Psi_{ij} = \Psi_i \circ g_{ij}$ and $\phi_{ij}(x) = 
\phi_i(x) Q_{ij}(\Psi_i^{-1}(x))\chi_{F_{ij}}(\Psi_i^{-1}(x))$, from $(2.5)$ and the assumption $(2.3)$,
for an appropriate rectangle $R_{ij}$ one gets
$$\chi_{D_{ij}^{k_i + 1}}(x) \phi_{ij}(\Psi_{ij}(x)) = \chi_{D_{ij}^{k_i+1}}(x) (\phi_i \circ \Psi_i \circ 
g_{ij}(x))Q_{ij}(g_{ij}(x))\chi_{F_{ij}}(g_{ij}(x)) $$
$$= \Phi_i(g_{ij}(x))Q_{ij}(g_{ij}(x))\chi_{F_{ij}}(g_{ij}(x))$$
$$= \chi_{R_{ij}}(x)\gamma_{ij}(x) \eqno (2.6)$$
This gives the version of part (2) of the main theorem that is needed in Theorem 2.1. Parts (3) and 
(4) are immediate, and we are done. 

\noindent {\bf Proof of Theorem 2.2:} The proof is by induction on the dimension $n$. When $n = 1$,
since 1-dimensional quasibump functions are smooth already the proof is straightfoward. Namely, one 
selects $r_1$ such that  $supp(\Psi) \subset [0,r_1]$. We then let
there be a single $F_j = \chi_{(0,r_1)}(x)$ and a single $Q_j(x) = \alpha(x)$ where 
$\alpha(x)= 1$ on $[0,r_1]$. The corresponding $g_j$ is just the identity map, and the case $n=1$ follows.
Assume now we that know Theorem 2.2 in dimension $n - 1$ and the hypotheses of
Theorem 2.2 hold for some $n$-dimensional situation. We break into two cases. 

\noindent {\bf Case 1:} The first case is when there are distinct monomials $t(x)$ and $u(x)$ and 
positive constants $c_1$ and $c_2$ such that whenever $\Phi(x) \neq 0$ we have
$$c_1  < {t(x) \over u(x)} < c_2 \eqno (2.7)$$
Write $t(x) = \prod_{i \in I} x_i^{l_i}$ and $u(x) = \prod_{i \in I'} x_i^{l_i}$, where each $l_i > 0$.
We can assume that $I \cap I' = \emptyset$ and $I \cup I' \neq \emptyset$. We change variables as 
follows. For a given $i \in I \cup I'$, let $m_i = \prod_{j \in I \cup I',\,j\neq i}l_j$ 
and let $x_i = y_i^{m_i}$.
If $i \notin I \cup I'$, let $x_i = y_i$. Let $x = g_1(y)$ be this coordinate change, and 
let $l$ denote $\prod_{i \in I \cup I'}l_i$. Then we have
$$t(g_1(y)) = \prod_{i \in I} y_i^l,\,\,\,u(g_1(y)) = \prod_{i \in I'} y_i^l \eqno (2.8)$$
Consequently, if $t_1(y) = \prod_{i \in I} y_i =  [t(g_1(y))]^{1 \over l}$, and $u_1(y) = 
\prod_{i \in I'} y_i = [u(g_1(y))]^{1 \over l}$, equation $(2.7)$ says that 
whenever $\Phi(g_1(y)) \neq 0$ one has
$$c_1^{1 \over l} < {t_1(y) \over u_1(y)} < c_2^{1 \over l} \eqno (2.9)$$
Effectively, we have reduced to the case when each $l_i$ is 1. We now prove Theorem 2.2 for case 1
by induction on $m = \min(|I|, |I'|)$. We start with the case where $m = 0$. In this case,
either $u_1(y)$ or $t_1(y)$ is a nonconstant monomial and by $(2.9)$ whenever $\Phi(g_1(y)) \neq 0$ 
that monomial is 
bounded below. Since the support of a quasibump function is also bounded, every $y_i$ appearing in this
monomial is therefore bounded below on the support $\Phi(g_1(y))$. Specifically, there is some 
$y_i$ and some constant $c$ such that $y_i > c$ whenever $\Phi(g_1(y))$ is nonzero. 

\noindent Since $\Phi(x)$ is a quasibump
function, so is $\Phi(g_1(y))$. As in $(2.1)$, we write $\Phi(g_1(y))$ as 
$$a(y) \prod_{l=1}^k b_l (c_l(y) {p_l(y) \over q_l(y)}) \eqno (2.10)$$
Because $y_i$ is bounded below, one can incorporate any powers of $y_i$ appearing in each $p_l(y)$ and 
$q_l(y)$ into the associated $c_l(y)$. Thus we may assume that the $p_l(y)$ and $q_l(y)$ do not depend on 
$y_i$. Let $c > 0$ be a constant such that $\prod_{l=1}^k b_l (c {p_l(y) \over q_l(y)}) = 1$ on the 
support of $\Phi(g_1(y))$. So we have
$$\Phi(g_1(y)) = \Phi(g_1(y))\prod_{l=1}^k b_l (c{p_l(y) \over q_l(y)})\eqno(2.11)$$
We apply the induction hypothesis to $\prod_{l=1}^k b_l (c{p_l(y) \over q_l(y)})$, a function of 
$n-1$ variables. Let $Q_j$ and $F_j$ be as in Theorem 2.2 applied to this function, and let $g_j^2(z)$
be the associated invertible monomial map, called $g_j(x)$ in the statement of that theorem. Then
$$\prod_{l=1}^k b_l (c{p_l(y) \over q_l(y)}) = \sum_j (\prod_{l=1}^k b_l (c{p_l(y) \over q_l(y)})
Q_j(y)\chi_{F_j}(y))\eqno (2.12)$$
Combining the last two equations gives
$$\Phi(g_1(y)) = \sum_j([\Phi(g_1(y))\prod_{l=1}^k b_l (c{p_l(y) \over q_l(y)})] Q_j(y)\chi_{F_j}(y))
\eqno (2.13)$$
Using this again on the bracketed expression, we have
$$\Phi(g_1(y)) = \sum_j\Phi(g_1(y)) Q_j(y)\chi_{F_j}(y)\eqno (2.14)$$
This will lead to our needed decomposition for $\Phi(g_1(y))$.
Next, note that we have
$$\Phi(g_1 \circ g_j^2(z)) = \sum_j \Phi(g_1(g_j^2(z)))Q_j(g_j^2(z))\chi_{F_j}(g_j^2(z)) \eqno (2.15)$$
Since Theorem 2.2 holds for $\prod_{l=1}^k b_l (c{p_l(z) \over q_l(z)})$ and $g_j^2(z)$, each 
${p_l(g_2^j(z)) \over q_l(g_2^j(z))}$ in $\Phi(g_1(g_j^2(z)))$ is a monomial. Similarly,
if as in Theorem 2.2 we write $Q_j(y) =\prod_l b_{jl} (c_{jl}(y) {p_{jl}(y) \over q_{jl}(y)})$,
we also have that each ${p_{jl}(g_2^j(z)) \over q_{jl}(g_2^j(z))}$ appearing in $Q_j(g_2^j(z))$ is a 
monomial. So since $\chi_{F_j}(g_2(z)) = 
\chi_{R_j}(z)$, equation $(2.15)$ shows that Theorem 2.2 holds for $\Phi$, with the associated 
invertible monomial maps given by $g_1 \circ g_2^j$. Hence we are done with the case where $m = 
min(|I|, |I'|) = 0$.

Now we assume $min(|I|, |I'|) = m > 0$, and that we know the theorem for $min(|I|, |I'|) = m - 1$.
For this fixed $m$, we induct on $max(|I|, |I'|)$ which we denote by $M$. The initial step $M = m$ is
done the same way as the inductive step, so we assume either $M = m$ or that $M > m$ and that we know
the result for $M - 1$. Let $i_1 \in I$ and let $i_2 \in I'$. Let $\alpha(t)$ be a $C^{\infty}$ 
function on $[0,\infty)$ that is equal to 1 for small enough $t$ and which satisfies $\alpha(t) + 
\alpha({1 \over t}) = 1$. In particular,
$$\alpha({y_{i_1} \over y_{i_2}}) + \alpha({y_{i_2} \over y_{i_1}}) = 1 \eqno (2.16)$$
Correspondingly, we decompose $\Phi(g_1(y)) = \Phi_1(y) + \Phi_2(y)$, where
$$\Phi_1(y) = \Phi(g_1(y))\alpha({y_{i_1} \over y_{i_2}}),\,\,\,\Phi_2(y) = \Phi(g_1(y))\alpha({y_{i_2} 
\over y_{i_1}}) \eqno (2.17)$$
Let $G_1(y)$ be the invertible monomial map whose $i_1$th component is $y_{i_1}y_{i_2}$, and whose $i$th
component is  $y_i$ for all $i \neq i_1$. Similarly, let $G_2(y)$ be the invertible monomial map whose 
$i_2$th component is $y_{i_1}y_{i_2}$, and whose $i$th component is  $y_i$ for all $i \neq i_2$. We have
$$\Phi_1(G_1(y)) = \Phi(g_1 \circ G_1(y))\alpha(y_{i_1}),\,\,\,\Phi_2(G_2(y)) = \Phi(g_1 \circ G_2(y))
\alpha(y_{i_2}) \eqno (2.18)$$
Both $\Phi_1(G_1(y))$ and $\Phi_2(G_2(y))$ are quasibump functions for the following reasons. First, since 
$\Phi(g_1(y))$ has compact support and for $k = 1,2$ the coordinate change $G_k$ is in the $y_{i_k}$ 
variable only, the $\Phi_k(G_k(y))$ have compact support in the $y_i$ variable for $i \neq i_k$. 
When $i = i_k$, the $\alpha(y_{i_k})$ factor 
ensures that $\Phi_k(G_k(y))$ has compact support in the $y_{i_k}$ variable as well. We conclude 
that both $\Phi_k(G_k(y))$ are compactly supported. Furthermore, since the $G_k$ are invertible 
monomial maps and $\Phi$ is a quasibump function, the $\Phi_k(G_k(y))$
are automatically of the form $(2.1)$. We conclude they are quasibump functions. 

Next, we translate
the condition $(2.9)$ in the new variables. Since $y_{i_1}$ appears in $t_1(y)$ and $y_{i_2}$
appears in $u_1(y)$, both of degree 1, the expression for ${t_1(G_1(y)) \over u_1(G_1(y))}$ is 
obtained from the expression for 
${t_1(y) \over u_1(y)}$ by removing  $y_{i_2}$ from the denominator, while ${t_1(G_2(y)) \over 
u_1(G_2(y))}$ is obtained from the expression for 
${t_1(y) \over u_1(y)}$ by removing  $y_{i_1}$ from the numerator, In either case, either $m = 
min(|I|, |I'|)$ or $M = max(|I|, |I'|)$ is reduced by 1. So by the inductive hypothesis
Theorem 2.2 applies to both $\Phi_1(G_1(y))$ and $\Phi_2(G_2(y))$. Let the decompositions
from Theorem 2.2 applied to these functions be given by
$$\Phi_1(G_1(y)) = \sum_j \Phi_1(G_1(y))Q_{ij}^1(y)\chi_{F_{ij}^1}(y),\,\,\,\,
\Phi_2(G_2(y)) = \sum_j \Phi_2(G_2(y))Q_{ij}^2(y)\chi_{F_{ij}^2}(y)\eqno (2.19)$$
Pulling back by $G_1$ and $G_2$ respectively, we have
$$\Phi_1(y) = \sum_j \Phi_1(y) Q_{ij}^1(G_1^{-1}y)\chi_{F_{ij}^1}(G_1^{-1}y),\,\,\,\,
\Phi_2(y) = \sum_j \Phi_2(y) Q_{ij}^2(G_2^{-1}y)\chi_{F_{ij}^2}(G_2^{-1}y)\eqno (2.20)$$
Since $\Phi(g_1(y)) = \Phi_1(y) + \Phi_2(y)$, we may add these, obtaining
$$\Phi(g_1(y)) = \sum_j \Phi_1(y)Q_{ij}^1(G_1^{-1}y)\chi_{F_{ij}^1}(G_1^{-1}y) + \sum_j \Phi_2(y)
Q_{ij}^2(G_2^{-1}y)\chi_{F_{ij}^2}(G_2^{-1}y) \eqno (2.21)$$
Going back to the original $x$ coordinates
$$\Phi(x) = \sum_j  \Phi_1(g_1^{-1}(x)) Q_{ij}^1((g_1 \circ G_1)^{-1}x) \chi_{F_{ij}^1}
((g_1 \circ G_1)^{-1}x) $$
$$+ \sum_j  \Phi_2(g_1^{-1}(x))
Q_{ij}^2((g_1 \circ G_2)^{-1}x) \chi_{F_{ij}^2}((g_1 \circ G_2)^{-1}x) \eqno (2.22)$$
Equation $(2.22)$ gives the required decomposition of the form $(2.2)$ for $\Phi(x)$. If $g_{ij}^1$, $g_{ij}^2$
denote the invertible monomial maps coming from Theorem 2.2 in $(2.19)$, then the invertible
monomial maps corresponding to the decomposition $(2.22)$ are given by
 $g_1 \circ G_1 \circ g_{ij}^1$ or $g_1 \circ G_2 \circ g_{ij}^2$. Since $\Phi_1(g_1^{-1}(x)) = 
\alpha({(g_1^{-1}(x))_{i_1} \over (g_1^{-1}(x))_{i_2}})\Phi(x)$ and $\Phi_2(g_1^{-1}(x)) = 
\alpha({(g_1^{-1}(x))_{i_2} \over (g_1^{-1}(x))_{i_1}})\Phi(x)$, 
the quasibump functions for $(2.22)$ (corresponding to the $Q_j(x)$ in $(2.2)$) are given by
$\alpha({(g_1^{-1}(x))_{i_1} \over (g_1^{-1}(x))_{i_2}})Q_{ij}^1((g_1 \circ G_1)^{-1}x)$ and
$\alpha({(g_1^{-1}(x))_{i_2} \over (g_1^{-1}(x))_{i_1}}) Q_{ij}^2((g_1 \circ G_2)^{-1}x)$, while the
sets corresponding to $F_j$ in $(2.2)$ are given by $(g_1 \circ G_1)^{-1}F_{ij}^1$ or $(g_1 \circ G_2)
^{-1}F_{ij}^2$. That $(2.22)$ satisfies the conclusions of Theorem 2.2 is a direct consequence of the
fact that the decomposition $(2.19)$ does. This completes the proof in Case 1.

\noindent {\bf Case 2:} The second case is when no expression of the form $(2.7)$ holds. Since
$\Phi(x)$ is a quasibump function, as in $(2.1)$ we can write
$$\Phi(x) = a(x) \prod_{l=1}^j b_l (c_l(x) {p_l(x) \over q_l(x)}) \eqno (2.23)$$
Since $\Phi(x)$ is compactly supported, we can multiply $(2.23)$ through by $\prod_{m=1}^n \alpha(x_m)$ 
for an appropriate function $\alpha(x)$ and not change the result. Hence we can assume
the $x_m$ are amongst the ${p_l(x) \over q_l(x)}$.
Write ${p_l(x) \over q_l(x)} = \prod_{i=1}^n x_i^{a_{il}}$. Here $a_{il}$ can be positive, negative, or
zero. By definition of quasibump function, there is a constant $c_1$ such that $\Phi(x)$ is supported 
on $\{x \in E: {p_l(x) \over
q_l(x)} < c_1$ for all $l \}$. If we define $(y_1,...,y_n)$ coordinates by $y_i = \ln(x_i)$, and
write the associated coordinate change as $x = e(y)$, then $\Phi(e(y))$ is supported on the set $A$ given by
$$A = \cap_{l=1}^j \{y \in \R^n: \sum_{i=1}^n a_{il} y_i < \ln(c_1) \} \eqno (2.24)$$
Define the set $A'$ by
$$A' = \cap_{l=1}^j  \{y \in \R^n: \sum_{i=1}^n a_{il} y_i < 0\} \eqno (2.25)$$
\noindent {\bf Lemma:} $A'$ is nonempty.

\noindent {\bf Proof:} We will show that if $A' = \emptyset$, then one must be in case 1 of this proof.  
 Suppose $A' = \emptyset$. Note that we may assume $\ln(c_1) > 0$; otherwise $A' \supset A \neq 
 \emptyset$. Since $A' = \emptyset$, we have $A = A - A'$ and therefore 
$$A \subset \cup_{l=1}^j  \{y \in \R^n: 0 < \sum_{i=1}^n a_{il} y_i < \ln(c_1)\} \eqno(2.26)$$
I claim that there is some $M> 0$ such that 
$$A \subset \{y \in \R^n:  -M < \sum_{i=1}^n a_{il} y_i < \ln(c_1)\} \eqno (2.27) $$ 
For suppose not. Then any $M > 0$ we can find a $y^l \in A$ such that $\sum_{i=1}^n a_{il} y_i^l < -M$. 
If $y$ denotes the average 
${1 \over j}\sum_{l=1}^jy^l$, then by the convexity of $A$, $y \in A$. Furthermore, since for all $l_1$ and $l_2$ 
we have $\sum_{i=1}^n a_{il_1} y_i^{l_2} < \ln(c_1)$, taking this average leads to
$$\sum_{i=1}^n a_{il} y_i = {1 \over j} \sum_{l=1}^j\sum_{i=1}^n a_{il} y_i^l < {j-1 \over j}\ln(c_1) - {M \over j}$$
Hence if $M$ is large enough, one has 
$$\sum_{i=1}^n a_{il} y_i < 0 \eqno (2.28)$$
Hence $y \in A'$, contradicting that $A' = \emptyset$. Therefore $(2.27)$ must hold. This automatically
implies that we are in case 1; for
the $y$ coordinates, $(2.7)$ translates to the existence of numbers $e_1$, $e_2$, $d_1,...,d_n$ such that
$e_1 < \sum_{i=1}^n d_i y_i < e_2$ whenever $(y_1,...,y_n) \in A$. Thus $(2.27)$ implies we in case 1 and 
the lemma is proved.

\noindent {\bf Lemma:} There is a vector $v$ such that $A \subset A' + v$

\noindent {\bf Proof:} let $w$ be any vector in $A'$. Then we have $\sum_{i=1}^n a_{il} w_i < 0$ for every 
$l$. If $y \in A$, for any $t > 0$, the vector  $y + tw$ satisfies 
$$\sum_{i=1}^n a_{il}(y + tw)_i = \sum_{i=1}^n a_{il} y_i + t (\sum_{i=1}^n a_{il} w_i) < \ln(c_1) + 
t (\sum_{i=1}^n a_{il} w_i)$$
Thus if $t$ is large enough, one has $\sum_{i=1}^n a_{il}(y + tw)_i < 0$ for all $y \in A$. In other words,
$y + tw \in A'$. Hence $A \subset A' - tw$ and we are done.

Since the coordinate functions $x_m$ are amongst the ${p_l(x) \over q_l(x)}$ (see below $(2.23)$),
the $n$ inequalities 
$\{y_m < 0\}$ are amongst the defining inequalities for $A'$ and therefore 
$$A' \subset \{y: y_m < 0 \hbox { for all } m\}\eqno (2.29)$$
We now "triangulate" $A'$. To be precise, by $(2.29)$ we see that $A' \cap \{y: \sum_{m=1}^n 
y_m = -1\}$ is a bounded convex 
polyhedron and therefore up to a set of measure zero can be written as a finite union 
$\cup_j S_j$ of $n-1$-dimensional simplices. For a given
$j$, we let $T_j$ be the union of all lines containing the origin and passing through $S_j$. We then have
$$A' =  \cup_j T_j $$
Since $A \subset A' + v$ for an appropriate vector $v$, we have 
$$A \subset \cup_j (T_j + v) \eqno (2.30)$$ 
Furthermore, each $T_j$ can be written in the form 
$$T_j = \cap_{l = 1}^n \{y \in \R^n: \sum_{i=1}^nh_{ilj} y_i < 0\}\eqno (2.31)$$
Equivalently,
$$T_j + v = \cap_{l = 1}^n \{y \in \R^n: \sum_{i=1}^nh_{ilj} y_i < \sum_{i=1}^nh_{ilj} v_i \}$$
We write $\eta_{lj} = \sum_{i=1}^nh_{ilj} v_i$, so that the above becomes
$$T_j + v = \cap_{l = 1}^n \{y \in \R^n: \sum_{i=1}^nh_{ilj} y_i < \eta_{lj}\}\eqno (2.32)$$
Next, we move everything back to the original $x$ coordinates. (Recall $y_i = \ln(x_i)$ for each $i$).
Denote the associated coordinate change by $x = e(y)$.
Writing $E = \{x: x_i > 0$ for all $i\}$,  we have 
$$e(T_j + v) = \cap_{l = 1}^n \{x \in E: \prod_{i=1}^n x_i^{h_{ilj}} < exp(\eta_{lj})\} \eqno (2.33)$$
By $(2.30)$, $e(A) \subset \cup_j e(T_j + v)$, while in $(2.24)$ we defined $A$ so that the original
quasibump function $\Phi(x)$ is nonzero only on $e(A)$. Hence if we denote $e(T_j + v)$ by $F_j$, we have
$$\Phi(x) = \sum_j \Phi(x)\chi_{F_j}(x)$$
For each $j$ let $Q_j(x)$ be some bump function of the form $\prod_{i=1}^n \alpha(x_i)$ that is 
equal to 1 on the support of $\Phi(x)$. One then has
$$\Phi(x) = \sum_j \Phi(x)Q_j(x)\chi_{F_j}(x) \eqno (2.34)$$
This will be the decomposition needed for Theorem 2.2. For a given $j$, the associated coordinate 
changes $g_j(z)$ are defined as follows. Let $H_j$ be the matrix whose $il$ entry is $h_{ilj}$,
and let $M_j = \{\epsilon_{dij}\}_{d,i=1}^n$ be a matrix such that  
integer coordinates such that $M_jH_j$ is $N$ times the identity matrix for a large integer $N$.
If one does the substitution $x_i = \prod_{d=1}^n z_d^{\epsilon_{dij}}$, then
$\prod_{i=1}^n x_i^{h_{ilj}}$ becomes $z_l^N$. We define $g_j(z)$ to be the map
$$ g_j(z) = (\prod_{d=1}^n z_d^{\epsilon_{d1j}},...,\prod_{d=1}^n z_d^{\epsilon_{dnj}})$$ 
Note that by $(2.33)$ we have
$$g_j^{-1}F_j = (0, exp({\eta_{1j} \over N})) \times .... \times (0, exp({\eta_{nj} \over N}))
\eqno (2.35)$$
\noindent {\bf Claim:} The $F_j$ and $g_j(z)$ defined this way satisfy the conclusions of Theorem 2.2.

\noindent {\bf Proof:} First, we check the $g_j(z)$ are invertible monomial maps; we have not
shown that each $\epsilon_{ijk}$ is nonnegative. Since $F_j \subset supp(\Phi)$, $F_j$ has compact
support. In particular each $x_i$ is bounded on $F_j$. Pulling back to the $z$ coordinates, this
means each $\prod_{d=1}^n z_d^{\epsilon_{dij}}$ is bounded on $(0, exp({\eta_{1j} \over N})) \times .... 
\times (0, exp({\eta_{nj} \over N}))$. This fact forces each $\epsilon_{dij}$ to be nonnegative. 
For on the curve $z = (t^{\alpha_1},...,t^{\alpha_n})$, $\alpha_d > 0$, the $i$th component of 
$g_j(z)$ is equal to $t^{\sum_d \alpha_d \epsilon_{dij}}$. If some $\epsilon_{dij}$ were
negative, by choosing $\alpha_d$ to be far larger than the other $\alpha_k$, the image of the curve 
under $g_j$ would go off to infinity.
We conclude that each $\epsilon_{dij}$ is nonnegative and therefore that
$g_j(z)$ is an invertible monomial map.

Next, observe that by $(2.34)$, part (1) of Theorem 2.2 holds. As for (2), we must show that 
if $p_k$ and $q_k$ are as in $(2.23)$, then each ${p_k(g_j(z)) \over q_k(g_j(z))}$ is a monomial.
Recall that by $(2.31)$, we have
$T_j = \cap_{l = 1}^n \{y \in \R^n: \sum_{i=1}^nh_{ilj} y_i < 0\}$, and that each equation 
${p_k(e(y)) \over q_k(e(y))} < 1$
is one of the defining equations for $A$. Translating the statement that  $T_j \subset A$ into the 
$x$ coordinates gives: For fixed $j$,
if $\prod_{i=1}^n x_i^{h_{ilj}} < 1$ for all $l$ then ${p_k(x) \over q_k(x)} < 1$ for all $k$. Pulling back by
$g_2$, we have that if $z \in (0, exp({\eta_{1j} \over N})) \times .... \times (0, exp({\eta_{nj} 
\over N}))$, then ${p_k(g_j(z)) \over q_k(g_j(z))} < 1$ for all $k$. Note that ${p_k(g_j(z)) \over q_k(g_j(z))}$ 
is of the form $\prod_{d=1}^n z_d^{\delta_{dkj}}$ for integers $\delta_{dkj}$. So exactly as in the last
paragraph, we must have that $\delta_{dkj} \geq 0$ for each $d$ and $k$. Hence for each $k$, 
${p_k(g_j(z)) \over 
q_k(g_j(z))}$ is a monomial. Since $j$ was arbitrary, we have part (2) of Theorem 2.2. Equation $(2.35)$
gives part (3) and we are done with the proof of Theorem 2.2.

\noindent {\bf 3. Proof of Theorem 1.1}

Most of this section will be devoted to proving Theorem 1.1. We will use the resolution of 
singularities algorithm
as described in Theorem 2.1 to reduce consideration to the simpler situation where there is only one
$f_i(x)$, and where that $f_i(x)$ is of the form $c(x)m(x)$ where $m(x)$ is a monomial and $c(x)$ is 
nonvanishing. Namely, we will use Theorem 2.1 to reduce things to proving the following:

\noindent {\bf Theorem 3.1:} Suppose $c(x)$ is a positive real-analytic function defined on a 
neighborhood of $[0,1]^n$, and $m(x) = \prod_{i=1}^n x_i^{m_i}$ is a nonconstant monomial. Let $Z$ 
denote $\{x \in [0,1]^n: m(x) = 0\}$. Define $K_{\phi}(t) $ by
$$K_{\phi}(t) = \int_{\{x \in (0,1)^n: c(x)m(x) < t\}} \phi(x)\,dx \eqno (3.1)$$
Then if $\phi(x)$ is a smooth function on $[0,1]^n$ supported in a sufficiently small neighborhood of 
$Z$, $K_{\phi}(t)$ has an asymptotic expansion of the following form, where the $\alpha$ 
range over an arithmetic progression of positive numbers:
$$K_{\phi}(t) = \sum_{\alpha \leq a}\sum_{i = 0}^{n-1} k_{i,\alpha}(\phi)\,t^{\alpha}\ln(t)^i
+ E_a(t) \eqno (3.2)$$
There are $M > 0$ and $D_{\alpha} > 0$ depending on $c(x)$ and $m(x)$ such that each $k_{i, \alpha}$ is a 
distribution with respect to $\phi$, supported on $Z$, satisfying
$$|k_{i,\alpha}(\phi)| < D_{\alpha} \sup_{|\beta| \leq M\alpha}\, \sup_{x \in Z} |\partial^{\beta}
\phi(x)| \eqno (3.3)$$
Also, there are $\epsilon > 0$ and $C_{a} > 0$ depending on $c(x)$ and $m(x)$ such that for 
any any $a$ and any $l$ satisfying $0 \leq l \leq a$ one has
$$|{d^l  \over dt^l}E_a(t)| < C_{a}\sup_{|\beta| \leq M(a + l)}\, \sup_{x \in (0,1)^n} |\partial^{\beta}
\phi(x)|\, t^{a + \epsilon - l} 
\eqno (3.4)$$

\noindent {\bf Proof that Theorem 3.1 implies Theorem 1.1:} Suppose we are in the setting of 
Theorem 1.1. We apply the resolution of singularities
algorithm simultaneously to each $f_l(x)$, each $g_l(x)$, and each difference $f_l(x) - f_m(x)$. Let
$U$ be as in the version of the main theorem of [Gr] from Theorem 2.1, and 
let $A = \{x \in \R^n: g_1(x) > 0,...,g_k(x) > 0\}$.  Let $\eta \in C_c(U)$ such that $\eta = 1$ on
a neighborhood $V$ of the origin, and let $\eta = \sum_i\eta_i$ be the decomposition given by the 
version of the main theorem given by Theorem 2.1. 

\noindent For any $\phi \in C_c(V)$, observe that $\phi = \sum_i \phi \eta_i$.
Thus we have
$$J_t = \int_{\{x \in A: 0 < f_1(x) < t,\,...\,,0 < f_l(x) < t\}} \phi(x)\,dx = 
 \sum_i \int_{\{x \in A: 0 < f_1(x) < t,\,...\,,0 < f_l(x) < t\}} \phi(x)\eta_i(x)\,dx$$
Let $\Psi_i$ be the composition of coordinate changes corresponding to $\eta_i$, let $B_i = 
\Psi_i^{-1}A$, and let $R_i$ be the rectangle $(0,a_1^i),...,(0,a_n^i)$
given by Theorem 2.1. Let $F_i = f_i \circ \Psi_i$. Then the $i$th term of the above expression 
becomes
$$\int_{\{x \in B_i \cap R_i: 0 < F_1(x) < t,\,...\,,0 < F_l(x) < t\}} \phi \circ \Psi_i(x) 
(\eta_i \circ \Psi_i(x))J_i(x)\,dx 
\eqno (3.5)$$
Here $J_i(x)$ denotes the Jacobian $\Psi_i$, and by Theorem 2.1 each $\eta_i \circ \Psi_i(x)$ is
smooth. Since each $g_j$ is being resolved by $\Psi_i(x)$, in
$(3.5)$ each $g_l \circ \Psi_i (x)$ is of the form $c_l(x)m_l(x)$, where $m_l(x)$ is a monomial and 
$c_l(x)$ is nonvanishing. In particular each $g_l \circ \Psi_i (x)$ is either everywhere positive or 
everywhere negative on $R_i$. Thus we have that either some $B_i \cap R_i = \emptyset$, whereupon
$(3.5)$ is zero, or that each $B_i\cap R_i = R_i$, whereupon $(3.5)$ is equal to
$$\int_{\{x \in R_i: 0 < F_1(x) < t,\,...\,,0 < F_l(x) < t\}} \phi \circ \Psi_i(x) 
(\eta_i \circ \Psi_i(x))J_i(x)\,dx \eqno (3.6)$$
Clearly we only have to consider terms of the form $(3.6)$. Next, notice that since we resolved all
differences $f_l(x) - f_m(x)$, each $F_l(x) - F_m(x)$ is also of the form $c_{lm}(x)m_{lm}(x)$, where 
$m_{lm}(x)$ is a monomial and $c_{lm}(x)$ is nonvanishing. As a result, each $F_l(x) - F_m(x)$ is
either everywhere positive or everywhere negative on $R_i$. Hence there is some $F_l(x)$ which is 
strictly larger than every other $F_m(x)$ everywhere on $R_i$, and $(3.6)$ is
$$\int_{\{x \in R_i: 0 < F_l(x) < t\}} \phi \circ \Psi_i(x) (\eta_i \circ \Psi_i(x))J_i(x)\,dx$$
Since the singularities of $g_l(x)$ are also being resolved, we can write $F_l(x) = c(x)m(x)$ and the 
above becomes
$$\int_{\{x \in R_i: 0 < c(x)m(x) < t\}} \phi \circ \Psi_i(x) (\eta_i \circ \Psi_i(x))J_i(x)\,dx \eqno (3.7)$$
To get $(3.7)$ into the form $(3.1)$, we do a scaling in $(3.7)$ to convert $R_i$ into
$(0,1)^n$. Denoting this scaling coordinate change by $(x_1,...,x_n) \rightarrow sx = 
(s_1x_1,...,s_nx_n)$, $(3.7)$ becomes
$$\prod_{l=1}^n s_l \int_{\{x \in (0,1)^n: 0 < c(sx)m(sx)) < t\}} \phi \circ \Psi_i(sx) 
(\eta_i \circ \Psi_i(sx))J_i(sx)\,dx$$
This integral is of the form $(3.1)$, with $\phi_i$ replaced by $\phi \circ \Psi_i(sx) 
(\eta_i \circ \Psi_i(sx))J_i(sx)$.
Hence applying Theorem 3.1 will give an asymptotic expansion for $(3.7)$ satisfying $(3.2) - (3.4)$
for $\phi \circ \Psi_i(sx) (\eta_i \circ \Psi_i(sx))J_i(sx)$ in place of $\phi_i(x)$. By the chain and 
product rules, $(3.2) - (3.4)$ then also
hold for $\phi$ itself. Adding over all $i$ gives Theorem 1.1 and we are done.

\noindent {\bf Proof of Theorem 3.1:} We induct on the number of positive $m_i$ in $m(x) = 
\prod_i x_i^{m_i}$. We start with the case where exactly
one $m_i$ is positive, and without generality we may assume $i = 1$, so that $c(x)m(x) = c(x)x_1^k$ for
some $k$. Since $c(x)$ is a positive real analytic function, we can write $c(x) = c_1(x)^k$ for some
positive real analytic $c_1$, and we have $c(x)m(x) = (c_1(x)x_1)^k$. We now do a coordinate change as
follows. For $i > 1$, we let $y_i = x_i$. For $i = 1$, we let $y_1 = c_1(x)x_1$. This is a smooth 
coordinate change on $x_1 < \delta$ for an appropriate $\delta$. Denote the coordinate change map by
$x = \Psi(y)$. Assuming $\phi(x)$ is supported on $x_1 < \delta$, we have
$$K_{\phi}(t) = \int_{\{x \in (0,1)^n: y_1 < t^{1 \over k}\}} \phi(\Psi(y))J(y)\,dy \eqno (3.8)$$
Here $J(y)$ denotes the Jacobian of this coordinate change. Thus $K_{\phi}(t)$
is the indefinite integral of $\Phi(y) = \phi(\Psi(y))J(y)$ in the $y_1$ variable from 0 to $t_1^{1/k}$. 
Hence if we substitute $\phi(\Psi(y))J(y) = \sum_{l = 0}^m  {\partial^l \Phi \over 
\partial y_1^l}(0,y_2,...,y_n){y_1^l \over l!} + O(y_1^{m+1})$ into $(3.8)$, one obtains 
$$K_{\phi}(t) = \sum_{l = 0}^m (\int_{(0,1)^{n-1}} {\partial^l \Phi \over \partial y_1^l}
(0,y_2,...,y_n)\,dy_2...\,dy_n) {t^{l+1 \over k} \over (l + 1)!} + O(t^{{m+2 \over k}})$$
This gives the desired asymptotic expression $(3.2)$ for
$K_{\phi}(t)$ in powers of $t^{1 \over k}$. The expressions $(3.3)$ and $(3.4)$ follow from the chain
rule applied to $\Phi(y) = \phi(\Psi(y))J(y)$. This concludes the proof of Theorem 3.1 for the case that
only one $m_i$ is positive.

Next, we assume that some $l > 1$ of the $m_i$ are positive, and we have shown the result for $l - 1$.
Without loss of generality, once again we assume that $m_1 = k > 0$. Like before we let $c_1(x)$ be
such that $c_1(x)^k$ $=c(x)$, so that $c(x)m(x) = (c_1(x)x_1)^k\prod_{i > 1}x_i^{m_i}$. Like above, 
for $x_1$ smaller than some $\delta$ we can do the coordinate change to variables $y_i$, where 
$y_1 = c_1(x)x_1$ and where $y_i = x_i$ for $i > 1$. Thus if $x = \Psi(y)$, we have
$$c(\Psi(y))m(\Psi(y)) = \prod_{i = 1}^n y_i^{m_i} \eqno (3.9)$$
Let $\xi(t)$ be a smooth function on $[0,\infty)$ that is supported in $[0, \delta)$, and which is equal
to 1 on a neighborhood of $0$. Correspondingly we write $\phi(x) = \phi_1(x) + \phi_2(x)$, where
$$\phi_1(x) = \phi(x)\xi(x_1),\,\,\,\phi_2(x) = \phi(x)(1 - \xi(x_1)) \eqno (3.10)$$
We correspondingly write $K_{\phi}(t) = K_{\phi}^1(t) + K_{\phi}^2(t)$, where
$$K_{\phi}^1(t) = \int_{\{x \in (0,1)^n: c(x)m(x) < t\}} \phi_1(x)\,dx,\,\,\, K_{\phi}^2(t) =
\int_{\{x \in (0,1)^n: c(x)m(x) < t\}} \phi_2(x)\,dx\eqno (3.11)$$
Note that $K_{\phi}^2(t)$ reduces to the case when $l - 1$ of the $m_i$ are positive, for the integrand
is supported on $x_1 > \delta '$ and thus the $x_1^k$ factor can be incorporated into the $c(x)$: One
does a linear variable change in the $x_1$ variable to turn $[\delta',1] \times [0,1]^{n-1}$ into
$[0,1]^n$ and then applies the $l - 1$ case. Hence it suffices to restrict our attention to 
$K_{\phi}^1(t)$. We change to the $y$ variables, obtaining
$$K_{\phi}^1(t) = \int_{\{x \in (0,1)^n: \prod_{i = 1}^n y_i^{m_i} < t\}} \phi_1(\Psi(y))J(y)
\,dy \eqno (3.12)$$
$J(y)$ is the Jacobian of the coordinate change. In $(3.12)$, the $c(x)$ factor has been removed from
the domain of integration, a fact that will make our arguments simpler. Write $\Phi(x) =
\phi_1(\Psi(y))J(y)$, and $(3.12)$ becomes
$$K_{\phi}^1(t) = \int_{\{x \in (0,1)^n: \prod_{i = 1}^n y_i^{m_i} < t\}} \Phi(y) \,dy \eqno (3.12')$$
\noindent {\bf Lemma 3.2:} To prove $K_{\phi}^1(t)$ has an asymptotic expansion satisfying the 
conditions of Theorem $(3.2)-(3.4)$, thereby proving Theorem 3.1, it suffices to find an asymptotic 
expansion for $K_{\phi}^1(t)$
satisfying the analogues of $(3.2)-(3.4)$ with $\phi$ replaced by $\Phi$. Furthermore, it suffices 
to consider only the case where $m_i > 0$ for all $i$.

\noindent {\bf Proof:} Since $\Phi(y) = \phi_1(\Psi(y))J(y) = \xi(\Psi_1(y))\phi(\Psi(y))J(y)$, by the 
chain and product rules, $|\partial^{\alpha}\Phi(x)|$ is bounded by 
$C_{\alpha}\sum_{|\beta| \leq |\alpha|}|\partial^{\beta}\phi(x)|$. Hence if the 
versions of $(3.2)-(3.4)$ hold for $K_{\phi}^1(t)$ with $\Phi$ in place of $\phi$, they will also hold 
for $K_{\phi}^1(t)$ with $\phi$ itself. This gives the first statement of Lemma 3.2.
As for the second statement, suppose
we have proved the asymptotics in the case where each $m_i > 0$. Then in the general situation,
we can fix those $y_i$ variables for which $m_i = 0$. Then the asymptotics hold for the 
integral in the remaining variables. By integrating the asymptotics with respect to the $y_i$ 
variables for which $m_i = 0$, one sees that the asymptotics hold for the original integral. 
This completes the proof of Lemma 3.2.

\noindent {\bf Lemma 3.3:} Suppose $g(y)$ is a $C^{\infty}$ function on $[0,1]^n$. Let $\sum_{\alpha}
g_{\alpha}y^{\alpha}$ denote the Taylor expansion for $g(y)$ about the origin. Then for each $N$ 
we can write 
$$g(y) = \sum_{\alpha_1,...,\alpha_n < N} g_{\alpha}y^{\alpha} + \sum_{\beta} 
{y^{\beta} \over \beta_1!..\beta_l!}[h_{\beta}(y_{l_{\beta}+1},...,y_n) - \sum_{i = 0}^{N-1} 
 {\partial^i h_{\beta} 
\over \partial y_{l_{\beta}+1}^i}(0,y_{l_{\beta}+2},...,y_n){y_{l_{\beta}+1}^i \over i!}] \eqno (3.13)$$
Here the sum in $\beta$ ranges over all $(\beta_1,...,\beta_l)$ with $0 \leq l < n$ and $0 \leq \beta_i 
< N$ for each $i$, $l_{\beta}$ denotes the number of entries $\beta$ has, and 
$h_{\beta}(y_{l_{\beta}+1},...,y_n)$ denotes $\partial^{\beta}g(0,..,0,y_{l_{\beta}+1},..,y_n)$. 

\noindent {\bf Proof:} Taylor expanding in $y_1$ we have
$$g(y) = \sum_{i = 0}^{N-1} {\partial^i g \over \partial y_1^i}(0,y_2,...,y_n) {y_1^i \over i!}+ 
(g(y) - \sum_{i = 0}^{N-1} {\partial^i g \over \partial y_1^i}(0,y_2,...,y_n) {y_1^i \over i!})
\eqno (3.14)$$
This gives the result for $n = 1$. When $n > 1$, one 
substitutes the $n-1$ dimensional case for ${\partial^i g \over \partial y_1^i}(0,y_2,...,y_n)$
into each term of the left series of $(3.14)$ and the lemma follows.

\noindent {\bf Lemma 3.4:} Let $p(y)$ denote one of the terms of the second sum of $(3.13)$; that is,
let $p(y)$ be of the form
$$p(y) = {y^{\beta} \over \beta_1!..\beta_l!} [h_{\beta}(y_{l_{\beta}+1},...,y_n) - \sum_{i = 0}^{N-1}
{\partial^i h_{\beta} 
\over \partial y_{l_{\beta}+1}^i}(0,y_{l_{\beta}+2},...,y_n){y_{l_{\beta}+1}^i \over i!}]$$
Let $j$ denote the index called $l_{\beta}+1$ in this term. Then if $\gamma$ is a multiindex 
such that $\gamma_j = 0$, then  
$$|\partial^{\gamma} p (y)| < C_{|\gamma|,N}'\,y_j^N||g||_{C^{nN + |\gamma|}}$$

\noindent {\bf Proof:} If one takes the $\gamma$ derivative of $p(y)$, one obtains some terms of the
 form 
$$cy^m[q(y_j,...,y_n) - \sum_{i = 0}^{N-1} {\partial^i q 
\over \partial y_j^i}(0,y_{j+1},...,y_n){y_j^i \over i!}] \eqno (3.15)$$
Here $y^m$ is a monomial, and $q$ is some partial derivative of the appropriate $h_{\beta}$ of order at most $|\gamma|$. 
Note that the bracketed expression is equal to $q(y_j,...,y_n)$ minus the first $N$ terms
of its Taylor expansion in the $y_j$ direction. As a result, by Taylor's theorem, the bracketed 
expression is equal to
${y_j^N \over N!} {\partial^N q  \over \partial y_j^N}(Y_j,y_{j+1},...,y_n)$ for some $Y_j$ between 
0 and $y_j$. Hence $(3.15)$ is bounded by $c'y_j^N||q||_{C^N}$. Since $q$ is a derivative of order at
most $|\gamma|$ of some $h_{\beta}$, which is itself a derivative of order at most $(n-1)N$ of $g$, we 
conclude that this term is bounded by $c''y_j^N||g||_{C^{nN + |\gamma|}}$. We add over all terms of 
$\partial^{\gamma} p(y)$ and we are done. 

\noindent We are now in a position to prove Theorem 1.1.

\noindent {\bf Proof of Theorem 1.1:} By lemma 3.2, it suffices to find an asymptotic expansion for $K_{\phi}^1(t)$
satisfying the analogue of $(3.2)-(3.4)$ with $\phi$ replaced by $\Phi$, and in $m(x) = \prod_i x_i^
{m_i}$ we may assume $m_i > 0$ for all $i$. Also, we can assume we know
the result inductively for dimensions $< n$. We fix some $a > 0$. With the goal of finding the
asymptotic expansion $(3.2)$ up to the power $t^a$, we apply Lemma 3.3 to $\Phi(y)$,
setting $N$ equal $\lfloor(a +  \epsilon)\max_i m_i \rfloor + 1$. We then
insert the result termwise into $(3.12')$. We obtain two types of terms. The first are terms of the form
$$ \Phi_{\alpha} \int_{\{y \in (0,1)^n: \prod_{i = 1}^n y_i^{m_i} < t\}} y^{\alpha} \,dy \eqno (3.16)$$
Here $\Phi_{\alpha}$ denotes the coefficient of $y^{\alpha}$ in $\Phi$'s Taylor expansion about the
origin. One can evaluate $(3.16)$ directly using calculus, by induction on the dimension for example, 
and prove that $(3.16)$ is of the form
$$\Phi_{\alpha} \sum_{\beta \leq M \alpha}\sum_{i = 0}^{n-1} c_{i,\beta}\,t^{\beta}\ln(t)^i 
\eqno (3.17)$$
This is the form that we need for our asymptotics. The second type of term we obtain comes from the
 second sum of $(3.13)$, a term denoted by $p(y)$ in Lemma 3.4. We write $p(y) = p_{y_j}(\bar{y})$,
where $y_j$ is as in the previous lemma and where $\bar{y}$ denotes the remaining $n-1$ variables. We
integrate the term first with respect to the
$\bar{y}$ variables and then with respect to the $y_j$ variable. The term becomes the following, where 
$\bar{m}(\bar{y})$ denotes the monomial ${m(y) \over y_j^{m_j}}$:
$$\int_0^1 \int_{\{y \in (0,1)^{n-1}: \bar{m}(\bar{y}) < {t \over y_j^{m_j}} \}} p_{y_j}(\bar{y})
\,d\bar{y}\,dy_j \eqno (3.18)$$
We next use the inductive hypothesis on the inner integral, and $(3.18)$ becomes
$$\int_0^1 [\sum_{\alpha \leq a}\sum_{k = 0}^{n-1} k_{k,\alpha}(p_{y_j}) ({t \over y_j^{m_j}})^{\alpha}
\ln({t \over y_j^{m_j}})^k + E_{a,y_j}({t \over y_j^{m_j}})] \,dy_j\eqno (3.19)$$
We expand $\ln({t \over y_j^{m_j}})^k = (\ln(t) - m_j \ln(y_j))^k$ and $(3.19)$ becomes
$$\int_0^1 [\sum_{\alpha \leq a}\,\,\sum_{k_1,k_2 = 0}^{n-1} \kappa_{k_1,k_2,\alpha}(p_{y_j}) 
({t \over y_j^{m_j}})^{\alpha}\ln(t)^{k_1}\ln(y_j)^{k_2} + E_{a,y_j}({t \over y_j^{m_j}})]\,dy_j \eqno (3.20)$$
Here each $\kappa_{k_1,k_2,\alpha}$ is a constant multiple of some $k_{k,\alpha}$. In 
particular, each $\kappa_{k_1,k_2,\alpha}$ satisfies $(3.3)$. We rewrite $(3.20)$ as 
$$\sum_{\alpha \leq a}\,\,\sum_{k_1,k_2 = 0}^{n-1} (\int_0^1 {\kappa_{k_1,k_2,\alpha}(p_{y_j})\ln(y_j)^{k_2}
\over y_j^{m_j\alpha}}\,dy_j) t^{\alpha}\ln(t)^{k_1}\,dy_j + \int_0^1 E_{a,y_j}({t \over y_j^{m_j}})\,dy_j
\eqno (3.21)$$
We will see that $(3.21)$ gives the desired asymptotic expansion. For by
$(3.3)$ and then Lemma 3.4 we have
$$|\kappa_{k_1,k_2,\alpha}(p_{y_j})| \leq  D_{\alpha} ||p_{y_j}||_{C^{M\alpha}(Z)}$$
$$\leq D_{\alpha}'y_j^N ||\Phi||_{C^{M\alpha + nN}(Z)} \eqno (3.22)$$
Hence, the absolute value of the $t^{\alpha}\ln(t)^{k_1}$ coefficient of $(3.21)$ satisfies 
$$|\sum_{k_2}\int_0^1 {\kappa_{k_1,k_2,\alpha}(p_{y_j})\ln(y_j)^{k_2} \over y_j^{m_j\alpha}}\,dy_j|\leq
D_{\alpha}'||\Phi||_{C^{M\alpha + nN}(Z)} \int_0^1 y_j^N   {|\ln(y_j)|^{n-1}\over y_j^{m_j\alpha}}
\,dy_j \eqno (3.23)$$
Since  $N \geq \max_k m_k(a +  \epsilon) \geq m_ja \geq m_j \alpha$, the right-hand integral above
is bounded, so $(3.23)$ is at most $C_{\alpha}'||\Phi||_{C^{M\alpha + nN}(Z)}$.  Since $\alpha \leq a$ and 
$N < M'a$ for some $M'$, we have $M\alpha + nN \leq (M + M'n)a$, and we conclude that $(3.23)$ is at most
$$C_{\alpha}' ||\Phi||_{C^{(M + M'n)a}(Z)}$$
$$\leq C_a''||\Phi||_{C^{(M + M'n)a}(Z)}$$
Here  $C_a''= \sup_{\alpha \leq a}C_{\alpha}'$. 
This gives us the desired estimate $(3.3)$ for the $t^{\alpha}\ln(t)^{k_1}$ coefficient.

We also have to analyze the error term $\int_0^1 E_{a,y_j}({t \over y_j^{m_j}})\,dy_j$. If $0 \leq l 
\leq a$, the $l$th derivative of this with respect to $t$ is equal to
$$\int_0^1 y_j^{-lm_j}{\partial^m E_{a,y_j} \over \partial t^m}({t \over y_j^{m_j}})\,dy_j \eqno (3.24)$$
Substituting $(3.4)$ into this, this is bounded in absolute value by
$$C_a \int_0^1 y_j^{-lm_j}||p_{y_j}||_{C^{Ma}((0,1)^{n-1})} ({t \over y_j^{m_j}})
^{a + \epsilon - l} \eqno (3.25)$$
By Lemma 3.4 this in turn is bounded by 
$$C_a\int_0^1 y_j^{-lm_j+ N}||\Phi||_{C^{Ma + nN}((0,1)^n)} ({t \over y_j^{m_j}})^{a + \epsilon - l}
\,dy_j$$
$$= C_a||\Phi||_{C^{Ma + nN}((0,1)^n)}(\int_0^1 y_j^{N - am_j - \epsilon m_j} 
 \,dy_j) t^{a + \epsilon - l}\eqno (3.26)$$
Since $N \geq am_j +\epsilon m_j$, $(3.26)$ is at most
$$C_a ||\Phi||_{C^{Ma + nN}((0,1)^n)}t^{a + \epsilon - l}\eqno (3.27)$$
And because $N \leq M'a$ for some $M'$, $(3.27)$ is bounded by 
$$C_a ||\Phi||_{C^{(M + M'n)a}((0,1)^n)}t^{a + \epsilon - l} \eqno (3.28)$$
This gives the desired estimate $(3.3)$ with $\Phi$ and we are done. To be clear, what we showed is 
that each $t^{\alpha}\ln(t)^j$'s coefficient is bounded in terms of $C^{M_0a}$
norms of $\Phi$ for $M_0 = M + M'n$, with corresponding estimates for the error terms. 
Note that for $(3.2)$-$(3.4)$ to hold we need
 $C^{M_0\alpha}$ norms, which makes a difference if $\alpha$ is a lot smaller than $a$. However,
if for a given $\alpha$ we consider the estimate obtained from the expansion to degree
$a = 2\alpha$ for example, we get the desired estimates.

{\noindent \bf Acknowledgements:} The author would like to thank D. H. Phong and B. Lichtin for helpful
comments.

{\noindent \bf References:}

\parskip=6pt

\noindent [AGV] V. Arnold, S Gusein-Zade, A Varchenko, {\it Singularities of differentiable maps
Volume II}, Birkhauser, Basel, 1988.

\noindent [At] M. Atiyah, {\it Resolution of singularities and division of distributions},
Comm. Pure Appl. Math. {\bf 23} (1970), 145-150.

\noindent [BGe] I. N Bernstein and S. I. Gelfand, {\it Meromorphy of the function ${\rm P}^
\lambda$}, Funkcional. Anal. i Prilo\v zen. {\bf 3} (1969), no. 1, 84-85.

\noindent [BiMi] E. Bierstone, P. Milman, {\it Semianalytic and subanalytic sets}, Inst. Hautes
Etudes Sci. Publ. Math. {\bf 67} (1988) 5-42.

\noindent [F] M.V. Fedoryuk, {\it The saddle-point method}, Nauka, Moscow, 1977.

\noindent [Gr] M. Greenblatt, {\it A Coordinate-dependent local resolution of singularities and 
applications}, to appear, J of Func. Analysis. 

\noindent [H1] H. Hironaka, {\it Resolution of singularities of an algebraic 
variety over a field of characteristic zero I},  Ann. of Math. (2) {\bf 79}
(1964), 109-203;

\noindent [H2] H. Hironaka, {\it Resolution of singularities of an algebraic 
variety over a field of characteristic zero II},  Ann. of Math. (2) {\bf 79}
(1964), 205-326. 

\noindent [J] P. Jeanquartier, {\it Développement asymptotique de la distribution de Dirac attachée 
à une fonction analytique}, (French) C. R. Acad. Sci. Paris Sér. A-B {\bf 201} (1970), 
A1159--A1161. 

\noindent [L] F. Loeser, {\it Volume de tubes autour de singularités}, (French) 
Duke Math. J. {\bf 53} (1986), no. 2, 443-455.

\noindent [M] B. Malgrange, {\it Integrales asymptotiques et monodromie}, Ann. Scient. Ecole Norm. 
Super., Ser. 4 {\bf 7} (1974) no. 3, 405-430. 

\noindent [PSSt] D. H. Phong, E. M. Stein, J. Sturm, {\it On the growth and 
stability of real-analytic functions}, Amer. J. Math. {\bf 121} (1999), no. 3, 519-554.

\noindent [S] E. Stein, {\it Harmonic analysis; real-variable methods,
orthogonality, and oscillatory integrals}, Princeton Mathematics Series Vol. 
43, Princeton University Press, Princeton, NJ, 1993.

\noindent [Va] V. Vassiliev, {\it The asymptotics of exponential integrals, Newton diagrams, and
classification of minima}, Functional Analysis and its Applications {\bf 11} (1977) 163-172.

\noindent Department of Mathematics, Statistics, and Computer Science \hfill \break
\noindent University of Illinois at Chicago \hfill \break
\noindent 851 S. Morgan Street \hfill \break
\noindent Chicago, IL 60607-7045 \hfill \break
\noindent greenbla@uic.edu
\end